\documentclass[graybox]{svmult}

    \usepackage{amsmath}
    \usepackage[english]{babel}
    \usepackage{amssymb}
    \usepackage{color}
    \usepackage{bigints}

\usepackage{mathptmx}       
\usepackage{helvet}         
\usepackage{courier}        
\usepackage{type1cm}        
\usepackage{makeidx}         
\usepackage{graphicx}        
\usepackage{multicol}        
\usepackage[bottom]{footmisc}

\makeindex             

\def \id{1\!\!\!1}

\def\graph#1{[\! [ #1]\! ]}

\def\wt{\widetilde}

\def\F{{\cal F}}

\def\G{{\cal G}}

\def\ff{{\mathbb F}}
\def\gg{{\mathbb G}}

\def\nn{{\mathbb N}}

\def\P{{\mathbb P}}

\def\E{\mathbb{E}}

\begin{document}
\title*{Projections, Pseudo-Stopping Times and the Immersion Property}
\titlerunning{Projections, pseudo-stopping times and the immersion property} 
\author{Anna Aksamit and Libo Li}
\institute{Anna Aksamit \at Mathematical Institute and Oxford Man Institute,
University of Oxford, 
United Kingdom, 
\email{anna.aksamit@maths.ox.ac.uk}
\and Libo Li \at Department of Mathematics and Statistics, 
University of New South Wales, 
Sydney, Australia, \email{libo.li@unsw.edu.au}}

\maketitle
\abstract*{
Given two filtrations $\ff\subset \gg$, we study when the $\ff$-optional projection and the $\ff$-dual optional projection coincide for a class of $\gg$-optional processes with integrable variation. It turns out that the later property is equivalent to the immersion property for $\ff$ and $\gg$ saying that
every $\ff$-local martingale is a $\gg$-local martingale, which, equivalently, may be characterised using the class of $\ff$-pseudo-stopping times. 
We also show that every $\gg$-stopping time can be decomposed into the minimum of two barrier hitting times.
}

\abstract{
Given two filtrations $\ff\subset \gg$, we study under which conditions the $\ff$-optional projection and the $\ff$-dual optional projection coincide for the class of $\gg$-optional processes with integrable variation.
It turns out that this property is equivalent to the immersion property for $\ff$ and $\gg$, that is
every $\ff$-local martingale is a $\gg$-local martingale, which, equivalently, may be characterised using the class of $\ff$-pseudo-stopping times. 
We also show that every $\gg$-stopping time can be decomposed into the minimum of two barrier hitting times.
}

\section{Introduction}

The study of pseudo-stopping times started in the paper by Williams \cite{W}. The author describes there an example of a non-stopping time $\tau$ which has the optional stopping property, namely, for every uniformly integrable martingale $M$, $\E(M_\tau)=\E(M_0)$.
Let us recall this example here.
Let $B$ be a Brownian motion and define:
$$T_1:=\inf \{t: B_t =1 \}\quad \textrm{and} \quad \sigma:=\sup\{t \leq T_1 : B_t =0\}.$$
Therefore $\sigma$ is the last zero of the process $B$ before it reaches one.
Let $\tau$ be the time of the maximum of $B$ over $[0, \sigma]$, that is
$$\tau:=\sup\{t<\sigma: B_t=B^*_t\}\quad \textrm{with} \quad
B^*_t:=\sup_{s\leq t}B_s.$$ 
Then, as shown in \cite{W}, $\tau$ has the optional stopping property. 
Such random times were then called pseudo-stopping times and further studied by Nikeghbali and Yor in \cite{NY1}. 
In particular, it was shown in \cite{NY1} that a finite random time $\tau$ is a pseudo-stopping time if and only if the optional projection of the process $\id_{[\![ \tau, \infty [\![}$ coincides with its dual optional projection. We want to study the conditions under which the later property holds not only for $\id_{[\![ \tau, \infty [\![}$ but for a larger class of processes. 
In other words, the main motivation of this work is to better understand the property that the optional projection is equal to the dual optional projection for processes of integrable variation, which is not true in general.

We work on a filtered probability space $(\Omega,\mathcal A,\mathbb F,\mathbb P)$, where $\mathbb{F}:=(\F_t)_{t\geq 0}$ denotes a filtration satisfying the usual conditions and we set $\F_\infty := \bigvee_{t\geq 0} \F_t \subset \mathcal{A}$. A process that is not necessarily adapted to the filtration $\ff$ is said to be {\it raw}. As convention, for any martingale, we work always with its c\`adl\`ag modification, while for any random process $(X_t)_{t\geq 0}$, we set $X_{0-} = 0$ and $X_\infty = \lim_{t\rightarrow \infty} X_t$ a.s, if it exists.

Let $\gg:=(\G_t)_{t\geq 0}$ be another filtration such that $\ff \subset \gg$, that is for each $t\geq 0$, $\F_t\subset \G_t$. 

The aim is to study under which conditions the $\ff$-optional projection and the $\ff$-dual optional projection coincide for the class of $\gg$-optional processes with integrable variation.
It turns out that this is closely related to the {\it immersion property} from the theory of enlargement of filtrations. 
A filtration $\ff$ is said to be {\it immersed} in $\gg$ and we write $\ff \hookrightarrow \gg$ if every $\ff$-local martingale is a $\gg$-local martingale. 
Often the immersion property is called the {\it hypothesis {\bf $({\cal H})$}} in the literature. 
We refer the reader to Br\'emaud and Yor \cite{BY} for further discussions and other conditions equivalent to the immersion property.

The results of this paper are also motivated by the study of the converse implications to the following known observations in the literature. 
Let us define a filtration $\ff^\tau:=(\F^\tau_t)_{t \geq 0}$ as the progressive enlargement of $\ff$ with $\tau$, i.e. the smallest right-continuous filtration containing $\ff$ such that $\tau$ is a stopping time, that is 
$$\F^\tau_t:=\bigcap_{s>t}(\F_s\lor \sigma(\tau\land s)).$$
In the reduced form approach to credit risk modelling (see Bielecki et al. \cite{BJR}), given a  filtration $\ff$, a popular way to model the default time $\tau$ is to use a barrier hitting time of an $\ff$-adapted increasing process and an independent barrier. 
It is known that a random time constructed in this fashion has the property that the filtration $\ff$ is immersed in $\ff^\tau$.
It is also known that the property that $\ff\hookrightarrow \ff^\tau$ implies that every $\ff^\tau$-stopping time is an $\ff$-pseudo stopping time. In fact, the authors of \cite{NY1} have observed that, given two filtrations $\ff$ and $\gg$ such that $\ff$ is immersed in $\gg$, every $\gg$-stopping time is an $\ff$-pseudo-stopping time.

The main contributions of this work are condition $(v)$ of Theorem \ref{ny2} and Theorem \ref{pseudoH}. 
In Theorem \ref{pseudoH} we show that the converse of the observation made by the authors of \cite{NY1} is true: if every $\gg$-stopping time is $\ff$-pseudo-stopping time then $\ff$ is immersed in $\gg$. 
Hence, it provides an alternative characterization of the immersion property based on pseudo-stopping times.
Furthermore in Theorem \ref{pseudoH} we show that the immersion property for $\ff$ and $\gg$ is equivalent to the property that the $\ff$-optional projection and $\ff$-dual optional projection coincide for the class of $\gg$-optional processes with integrable variation. 

As an application of Theorem \ref{ny2} $(v)$, 
which gives another equivalent characterisation for pseudo-stopping times, and Theorem \ref{pseudoH}, we provide in Proposition \ref{honestnotpseudo} an alternative proof  to a result regarding the immersion property and the progressive enlargement with honest times. The advantage of our method is that we do not use specific structures of the progressive enlargement and the characterization of predictable sets as done in Jeulin \cite{J2}.

As another application, assuming that $\ff$ is immersed in $\gg$, we show in Theorem \ref{gstoping} that every $\gg$-stopping time can be written as the minimum of two $\gg$-stopping times, 
one of which is a barrier hitting time of an $\ff$-adapted increasing process, where the barrier is {\it 'almost'} independent, and the other is an $\ff$-pseudo-stopping time whose graph is contained in the union of the graphs of a family of $\ff$-stopping times. 

\section{Characterisation of pseudo-stopping times}

The main object of interest in this section is the class of pseudo-stopping times. We start with recalling the definition of pseudo-stopping times from \cite{NY1}, with a slight modification, that is the random time is allowed to take the value infinity.

\begin{definition}\label{pseudodef}
A random time $\tau$ is an {\it $\ff$-pseudo-stopping time} if for every uniformly integrable $\ff$-martingale $M$, we have $\E(M_\tau)= \E(M_0)$.
\end{definition}

The main tools used in this study are the (dual) optional projections onto the filtration $\ff$. 
We record here some known results from the general theory of stochastic processes. 
For more details on the theory the reader is referred to He et al. \cite{HWY} or Jacod and Shiryaev \cite{JS} and for specific results from the theory of enlargement of filtrations to Jeulin \cite{J2}. 

For any locally integrable variation process $V$, we denote the $\ff$-optional projection of $V$ by $^{o}V$ and the $\ff$-dual optional projection of $V$ by $V^{o}$.
It is known that the process $N^V:= \,^oV - V^o$ is a uniformly integrable $\ff$-martingale with $N^V_0 = 0$ and $\,^o(\Delta V) = \Delta V^o$. 

We specialize the above notions to the study of random times. 
For an arbitrary random time $\tau$, we set $A:=\id_{[\![\tau,\infty[\![}$
and define \hfill\break 
\noindent $\bullet$ the supermartingale $Z$ associated with $\tau$, $Z :=\,^{o}(\id_{[\![0,\tau[\![})= 1- \,^{o}A$,\hfill\break 
\noindent $\bullet$ the supermartingale $\widetilde Z$ associated with $\tau$, $\wt Z :=\,^{o}(\id_{[\![0,\tau]\!]})= 1- \,^{o}(A_-)$,\hfill\break 
\noindent $\bullet$ the martingale $m :=1-\left ( \,^{o}A-A^o\right )$.

These processes are linked through the following relationships:
\begin{equation*}
\label{relation}
Z=m-A^o \quad \textrm{and} \quad \widetilde Z=m- A^o_-.
\end{equation*}

We present in Theorem \ref{ny2} an extension of Theorem 1 from Nikeghbali and Yor \cite{NY1}. We extend their result in two directions.
Firstly, we allow for non-finite pseudo-stopping times. 
Secondly, Theorem 1 from \cite{NY1} states that if either all $\ff$-martingales are continuous or the random time $\tau$ avoids all {\rm finite} $\ff$-stopping times, i.e. $\Delta A^o  = 0$, then the random time $\tau$ is an $\ff$-pseudo-stopping time if and only if the process $Z$ is a decreasing $\ff$-predictable process.
We will remove these additional assumptions and present another equivalent characterization based on the process $\wt Z$ instead of $Z$ in condition $(v)$ of Theorem \ref{ny2}. 
We point out that the equivalence of condition $(v)$ of Theorem \ref{ny2} is one of the key results of this paper.

Before presenting Theorem \ref{ny2} let us give a motivating example of a random time which is not a pseudo-stopping time and $Z= \wt Z$ is decreasing but not predictable. 
It illustrates the importance of the {\it c\`agl\`ad} property in condition $(v)$ of Theorem \ref{ny2}.
\begin{example}
\label{ex}
Let $N$ be a Poisson process with intensity $\lambda$ and jump times $(T_n)_n$. 
Consider the random time $\tau=\frac{1}{2}(T_1+T_2)$. 
Then we obtain
\begin{equation*}
\E(N_{\tau\land 1}-\lambda (\tau\land 1))<\E(N_{T_1\land 1}-\lambda(T_1\land 1))=0
\end{equation*}
which implies that $\tau$ is not a pseudo-stopping time.
Furthermore we compute
$$\wt  Z_t=Z_t
  =  \id_{\{T_1>t\}}+\id_{\{T_1\leq t\}} \id_{\{T_2>t\}} e^{-\lambda  (t-T_1)}
$$
hence $\wt Z=Z$ is a decreasing and c\`adl\`ag process. 
This example is further studied in Proposition 5.3 in \cite{ACDJ}.
\end{example}

\begin{theorem}
\label{ny2}
The following are equivalent:\\
(i) $\tau$ is an $\ff$-pseudo-stopping time;\\
(ii) $A^o_\infty=\P(\tau<\infty\,|\,\F_\infty)$;\\
(iii) $m=1$ or equivalently $\,^oA= A^o$;\\
(iv) for every $\ff$-local martingale $M$, the process $M^\tau$ is an $\ff^\tau$-local martingale;\\
(v) the process $\wt Z$ is a c\`agl\`ad decreasing process.
\end{theorem}

Before proceeding to the proof of Theorem \ref{ny2}, we give an auxiliary lemma which characterizes the main property of our interest, that is, given a process of finite variation, when is its optional projection equal to the dual optional projection.
\begin{lemma}\label{Hloc}
Given a raw locally integrable increasing process $V$, the following are equivalent:\\
(i) $\,^o(V_-)$ is a c\`agl\`ad increasing process;\\
(ii) $\,^o(V_-) = V^o_-$; \\
(iii) $\,^o(V_-) = \,^o V_-$; \\
(iv) $\,^o V = V^o$ .
\end{lemma}

\begin{proof}
For any raw locally integrable increasing process $V$, from classic theory we know that the process $N^V:= \,^o V - V^o$ is a uniformly integrable martingale with $N^V_0 = 0$ and $\,^o(\Delta V) = \Delta V^o$. 
As a consequence we have 
\begin{align}
N^V= \,^o(V_-) - V^o_- \quad \mathrm{and} \quad N^V_-= \,^o V_- - V^o_- .\label{m1}
\end{align}
If $\,^o(V_-)$ is a c\`agl\`ad increasing process, then from \eqref{m1}, we see that $N^V$ is a predictable martingale of finite variation, therefore is constant and equal to zero, since predictable martingales are continuous which shows $(i) \implies (ii)$ and $(i) \implies (iv)$.
To prove $(iv) \implies (ii)$, it is enough to use the definition of $N^V$. Since $N^V$ is c\`adl\`ag, we know that $N^V\equiv 0$ if and only if $N^V_-\equiv 0$. This fact combined with \eqref{m1} gives the equivalence between $(ii)$ and $(iii)$ and the equivalence between $(ii)$ and $(iv)$.
\qed \end{proof}

\noindent \emph{Proof of Theorem \ref{ny2}.}
To see that $(i)$ and $(ii)$ are equivalent, suppose $\tau$ is an $\ff$-pseudo-stopping time. 
Then, by properties of optional and dual optional projection, for any uniformly integrable $\ff$-martingale $M$ we have
\[
\E(M_\tau \id_{\{\tau<\infty\}})= \E(\int_{[0,\infty)} M_{s}dA^o_s) = \E(M_\infty A^o_\infty ).
\]
Therefore, the equality, $\E(M_\tau)=\E(M_\infty)$ holds true for every uniformly integrable $\ff$-martingale $M$ if and only if $A^o_\infty=\P(\tau<\infty\,|\,\F_\infty)$, since
$\E(M_\tau)=\E(M_\infty(A^o_\infty+ \P(\tau=\infty|\F_\infty)))$. 

On the other hand, we have $\,^oA_\infty = \lim_{s\rightarrow\infty} \P(\tau\leq s\,|\,\F_s) = \P(\tau < \infty\,|\,\F_\infty)$ a.s., and from the definition of $m$, we note that $(ii)$ holds if and only if $(iii)$ holds, that is $m = 1$ or equivalently $\,^oA= A^o$. The equivalence of $(iii)$ and $(v)$ follows directly from Lemma \ref{Hloc}.

To see that $(i)\implies (iv)$, let $M$ be a uniformly integrable $\ff$-martingale.  
For any $\ff^\tau$-stopping time $\nu$, from Dellacherie et al. \cite{DMM} (page 186), we know there exists an $\ff$-stopping time $\sigma$ such that $\tau\wedge \nu = \tau\wedge \sigma$. 
Therefore, from the definition of pseudo-stopping time,
\begin{equation*}
\E(M_{\tau\wedge \nu}) = \E (M_{\tau\wedge \sigma}) = \E(M_0),
\end{equation*}
which shows that $M^\tau$ is a uniformly integrable $\ff^\tau$-martingale by Theorem 1.42 \cite{JS}.
The implication $(iv)\implies (i)$ is straightforward.
\qed 

\begin{remark}
The importance of the c\`agl\`ad property in condition $(v)$ of Theorem \ref{ny2} is illustrated in Example \ref{ex}. 
From this example we also see that a decreasing supermartingale $Z$ is not sufficient to ensure that the time is a pseudo-stopping time.
We would also like to point out that condition $(v)$ in Theorem \ref{ny2} is crucial when working with non-continuous filtrations and it is used later in the proof of Proposition \ref{honestnotpseudo}.
\end{remark}

\section{Main results and applications}

In this section, we formulate in Theorem \ref{pseudoH} our main result which provides the necessary and sufficient conditions for the property that the $\ff$-dual optional projection and $\ff$-optional projection of any $\gg$-optional process of integrable variation coincide. 
As a part of this result we derive a new characterization of the immersion property in terms of pseudo-stopping times. 

\begin{theorem}\label{pseudoH}
Given filtrations $\ff$ and $\gg$ such that $\ff \subset \gg$, the following are equivalent,\\
(i) the $\ff$-dual optional projection of any $\gg$-optional process of integrable variation is equal to its $\ff$-optional projection;\hfill\break
(ii) every $\gg$-stopping time is an $\ff$-pseudo-stopping time; \hfill\break
(iii) the filtration $\ff$ is immersed in $\gg$.
\end{theorem}

\begin{proof}
The implication $(i) \implies (ii)$ follows directly from Theorem \ref{ny2} $(iii)$. 

Let us now show the implication $(iii) \implies (i)$. 
Under the immersion property, the $\ff$-optional projection of any bounded $\gg$-optional process is equal to its optional projection on to the constant filtration $\F_\infty$ (see Bremaud and Yor \cite{BY}). 
More explicitly, for any given locally integrable increasing $\gg$-adapted process $V$, we have $\,^{o}(V_{-})_\sigma = \E(V_{\sigma-}|\F_\infty)$ for any $\ff$-stopping time $\sigma$. 
From this we see that the process $\,^{o}(V_-)$ is increasing c\`agl\`ad and $(i)$ follows from Lemma \ref{Hloc}.

To show $(ii) \implies (iii)$, suppose that $M$ is a uniformly integrable $\ff$-martingale and $\nu$ is any $\gg$-stopping time.
Since every $\gg$-stopping time is an $\ff$-pseudo-stopping time, we have $\E(M_\nu) = \E(M_0)$ for every $\gg$-stopping time $\nu$, which by Theorem 1.42 in \cite{JS}, implies that $M$ is a uniformly integrable $\gg$-martingale.

The theorem is now proved, however, for the sake of completeness, let us directly show that $(iii) \implies (ii)$. To this end, let $M$ be any uniformly integrable $\ff$-martingale and $\nu$ a $\gg$-stopping time. Then, from the immersion property, $M$ is a uniformly integrable $\gg$-martingale and $\E(M_\nu) = \E(M_0)$, which implies $\nu$ is an $\ff$-pseudo-stopping time. 
\qed \end{proof}

We now give two applications of our main results in Theorem \ref{ny2} and Theorem \ref{pseudoH}. 
An important class of random times is the class of honest times. A random time $\tau$ is an $\ff$-honest time if for every $t > 0$ there exists an $\F_t$-measurable random variable $\tau_t$ such that $\tau=\tau_t$ on $\{\tau<t\}$. In Proposition \ref{honestnotpseudo} we relate pseudo-stopping times with honest times and, as an application of Theorem \ref{ny2} $(v)$ combined with the equivalence between $(ii)$ and $(iii)$ in Theorem \ref{pseudoH}, we recover a new proof of a result regarding honest times and the immersion property found in Jeulin \cite{J2}. 
Therein the result is obtained by computing explicitly the $\gg$-semimartingale decompositions of $\ff$-martingales.
The equivalence $(i) \iff (ii)$ in Proposition \ref{honestnotpseudo} was already presented in Proposition 6 in \cite{NY1} under the simplifying assumption that all $\ff$-martingales are continuous and the proof therein uses distributional arguments. 
Here, we show that a similar result can be obtained in full generality by using sample path properties based on Theorem \ref{ny2} $(v)$.
\begin{proposition}
\label{honestnotpseudo}
Let $\tau$ be a random time. The following conditions are equivalent,\\
(i)  $\tau$ is equal to an $\ff$-stopping time on $\{\tau<\infty\}$,\\
(ii) $\tau$ is an $\ff$-pseudo-stopping time and an $\ff$-honest time.\\
In particular if $\tau$ is an $\ff$-honest time which is not equal to an $\ff$-stopping time on $\{\tau<\infty\}$ and a $\gg$-stopping time for some filtration $\gg\supset \ff$ then $\ff$ is not immersed in $\gg$. 
\end{proposition}
\begin{proof}
The implication $(i) \Longrightarrow (ii)$ is obvious so we show only $(ii) \Longrightarrow (i)$. Given that $\tau$ is a honest time, by Proposition 5.2. in \cite{J2}, we have that $\tau = \sup\{t: \wt Z_t=1 \}$ on $\{\tau<\infty\}$. On the other hand, by Theorem \ref{ny2} $(v)$, the pseudo-stopping time property of $\tau$ implies that $\wt Z=1-A^o_-$. Therefore, on $\{\tau<\infty\}$,
\begin{align*}
\tau&=\sup\{t: \wt Z_t=1 \} =\sup\{t: A^o_{t-}=0 \} =\inf\{t:A^o_{t}>0 \},
\end{align*}
so, $\tau$ is equal to an $\ff$-stopping time on $\{\tau<\infty\}$.\\
Therefore if $\tau$ is an $\ff$-honest time which is not equal to an $\ff$-stopping time on $\{\tau<\infty\}$ and a $\gg$-stopping time for some filtration $\gg\supset \ff$ then, by Theorem \ref{pseudoH}, $\ff$ is not immersed in $\gg$. 
\qed \end{proof}

In the remaining, given that $\ff \hookrightarrow \gg$, we show that every $\gg$-stopping time can be written as the minimum of two barrier hitting times for which the $\F_\infty$-conditional distribution of the barriers can be computed.
The proof of our final result given in Theorem \ref{gstoping} relies on the equivalence $(i) \iff (iii)$ in Theorem \ref{pseudoH}.

\begin{theorem}\label{gstoping}
Assume that $\ff\hookrightarrow \gg$ and let $\tau$ be a $\gg$-stopping time. 
Then $\tau$ can be written as $\tau_c\wedge\tau_d$, where:

\noindent (i) The random time $\tau_c$ is a $\gg$-stopping time which avoids all finite $\ff$-stopping times.
Denote by $A^{c,o}$ the $\ff$-dual optional projection of the process $\id_{[\![\tau_c,\infty[\![}$.
Then the $\F_\infty$-conditional distribution of $A^{c, o}_{\tau_c}$ is uniform on the interval $[0,A^{c, o}_\infty)$, with an atom of size $1-A^{c, o}_\infty$ at $A^{c, o}_\infty$, 
that is
\[
\P(A^{c, o}_{\tau_c}\leq u|\F_\infty)  
=  u\id_{\{u < A^{c, o}_\infty\}} + \id_{\{u \geq A^{c, o}_\infty\}}.
\]
\hfill\break
\noindent (ii) The random time $\tau_d$ is a $\gg$-stopping time whose graph is contained in the disjoint union of the graphs of the jump times of the process $A^o$ given by $(\sigma_k)_{k\in\nn}$.
Denote by $A^{d,o}$ the $\ff$-dual optional projection of the process $\id_{[\![\tau_d,\infty[\![}$.
Then
\begin{align*}
\P(A^{d, o}_{\nu_d}= u|\F_\infty)
&= \sum_k \id_{\{A^{d, o}_{\sigma_k}= u\}}\Delta A^{d,o}_{\sigma_k}.
\end{align*}
\end{theorem}

Before proceeding to the proof of Theorem \ref{gstoping}, we show that in fact any random time $\tau$ can be written as a barrier hitting time of an $\ff$-adapted increasing process given the appropriate barrier. 
We refer the reader to Remark 3.2 in Gapeev \cite{G} where the author considers the situation where the process $A^o$ is strictly increasing. 
We will demonstrate this result with no assumptions on $A^o$.
\begin{lemma}
A random time $\tau$ can be written as the barrier hitting time of the process $A^o$ with the barrier $A^o_\tau$, that is $\tau = \inf\{t>0 : A^o_t \geq A^o_\tau\}.$
\end{lemma} 
\begin{proof}
We first define another random time $\tau^*$ by setting
\[
\tau^* := \inf\{t>0 : A^o_t \geq A^o_\tau\}.
\]
To see that $\tau^* = \tau$ (it is obvious that $\tau^* \leq \tau$), we use Lemma 4.2 of \cite{J2} which states that the left-support of the measure $dA$, i.e.,  
\[
\{(\omega, t) : \forall \varepsilon >0 \quad A_t(\omega)>A_{t-\varepsilon}(\omega)\}=\graph{\tau}
\]
 belongs to the left-support of $dA^o$, i.e., to the set $\{(\omega, t) : \forall \varepsilon >0 \quad A^o_t(\omega)>A^o_{t-\varepsilon}(\omega)\}$.
 \qed \end{proof}

\noindent \emph{Proof of Theorem \ref{gstoping}.}
For any $\gg$-stopping time $\tau$ and the set $D: ={\{\Delta A^o_\tau>0\}}\in \G_\tau$, we see that $\tau$ can be written as $\tau_c\wedge \tau_d$, where $\tau_c:= \tau\id_{D^c}+\infty\id_{D}$ and $\tau_d:= \tau\id_{D}+\infty\id_{D^c}$. The random times $\tau_c$ and $\tau_d$ are therefore $\gg$-stopping times, where $\tau_c$ avoids finite $\ff$-stopping times and the graph of $\tau_d$ is contained in the graphs of the jump times of $A^o$. For more details on this decomposition of a random time see \cite{ACJ}.

Given $\tau$ is a $\gg$-stopping time that avoids all finite $\ff$-stopping times. The $\F_\infty$-conditional distribution of $A^o_\tau$ is given by
\[
\E(\id_{\{A^o_\tau\leq u\}}|\F_\infty)  = 
\E(\id_{\{A^o_\tau\leq u\}}|\F_\infty) \id_{\{u < A^o_\infty\}}  + \id_{\{u \geq A^o_\infty\}} .
\]
Let us set $C$ to be the right inverse of $A^o$, then the first term in the right hand side above is
\begin{align*}
\E(\id_{\{A^o_\tau\leq u\}}\id_{\{C_u < \infty\}}|\F_\infty) & = \E(\id_{\{\tau\leq C_u\}}\id_{\{C_u < \infty\} }|\F_{C_u})\\
                    & = \,^{o}A_{C_u}\id_{\{C_u < \infty\}}\\
		    & = A^o_{C_u}\id_{\{C_u < \infty\}}\\
                    & = u\id_{\{u < A^o_\infty\}}
\end{align*}
where we apply Theorem \ref{pseudoH} in the third equality, while the last equality follows from the fact that $A^o_{C_u} = u$, since $A^o$ is continuous except, perhaps, at infinity. This implies that the $\F_\infty$-conditional distribution of $A^{o}_\tau$ is uniform on $[0,A^o_\infty)$.

On the other hand, given $\tau$ is a $\gg$-stopping time whose graph is contained in the graphs of the jump times of $A^o$ given by $(\sigma_k)_{k\in \nn}$. Then
\begin{align*}
\P(A^o_\tau= u|\F_\infty) &= 
\sum_k \P(\{\tau = \sigma_{k}\}\cap \{A^o_{\sigma_k}= u\}|\F_\infty)\\
				       &= \sum_k \id_{\{A^o_{\sigma_k}= u\}}\P(\tau=\sigma_{k}|\F_\infty)\\
				&= \sum_k \id_{\{A^o_{\sigma_k}= u\}}\Delta A^o_{\sigma_k}       
\end{align*}
where the last equality follows from the fact that $\ff\hookrightarrow \gg$.
\qed 

\begin{remark}
As a special case of Theorem \ref{gstoping}, if $\tau$ is a finite $\gg$-stopping time that avoids finite $\ff$-stopping times, then $A^o_\tau$ is independent of $\F_\infty$ and uniformly distributed on the interval $[0,1]$. 
In this case, the $\gg$-stopping time $\tau$ is a barrier hitting time of an $\ff$-adapted increasing process, with the barrier being independent from $\F_\infty$. 
This is a class of random times widely used in credit risk modelling to model default times.
\end{remark}

\noindent {\bf Acknowledgement:} The authors wish to thank Monique Jeanblanc, Jan Ob\l\'{o}j, Marek Rutkowski and the anonymous referee for their careful readings and valuable advices on the writing of this paper and acknowledge the generous financial supports of {\it 'Chaire March\'es en Mutation'}, F\'ed\'eration Bancaire Fran\c caise. 
The second author also wishes to thank Arturo Kohatsu-Higa for his hospitality in Japan and acknowledge the grants of the Japanese government. 
The research of the first author has been partially supported by the European Research Council under the European Union's Seventh Framework Programme (FP7/2007-2013) / ERC grant agreement no. 335421.

\end{document}